\documentclass[leqno,english]{smfart}

\usepackage[leqno]{amsmath}
\usepackage{amssymb, amsthm, enumerate, latexsym, mathrsfs, mdwlist, stmaryrd}
\usepackage{ae,aecompl}
\usepackage[all]{xy}

\theoremstyle{remark}
\newtheorem{nul}{}[section]

\newtheorem*{rem*}{Remark}

\theoremstyle{definition}
\newtheorem{dfn}[nul]{Definition}

\newtheorem{ntn}[nul]{Notation}
\newtheorem{exm}[nul]{Example}

\newtheorem*{dfn*}{Definition}
\newtheorem*{axm*}{Axiom}
\newtheorem*{ntn*}{Notation}
\newtheorem*{exm*}{Example}
\newtheorem*{exr*}{Exercise}
\newtheorem*{int*}{Intuition}
\newtheorem*{qst*}{Question}

\theoremstyle{plain}

\newtheorem{thm}[nul]{Theorem}
\newtheorem{prp}[nul]{Proposition}
\newtheorem{cor}[nul]{Corollary}
\newtheorem{lem}[nul]{Lemma}

\newtheorem*{thm*}{Theorem}
\newtheorem*{prp*}{Proposition}
\newtheorem*{cor*}{Corollary}
\newtheorem*{lem*}{Lemma}
\newtheorem*{cnj*}{Conjecture}

\numberwithin{equation}{nul}

\swapnumbers

\DeclareMathOperator{\Aut}{Aut}

\DeclareMathOperator{\const}{const}
\DeclareMathOperator{\colim}{colim}

\DeclareMathOperator{\End}{End}
\DeclareMathOperator{\END}{\underline{\End}}

\DeclareMathOperator{\mor}{mor}
\DeclareMathOperator{\Mor}{Mor}
\DeclareMathOperator{\MOR}{\underline{\Mor}}

\DeclareMathOperator{\Obj}{Obj}

\newcommand{\coprd}{\amalg}

\newcommand{\MM}{\mathbf{M}}
\newcommand{\NN}{\mathbf{N}}

\newcommand{\VV}{\mathbf{V}}

\newcommand{\XX}{\mathbf{X}}

\newcommand{\op}{\mathrm{op}}

\usepackage[naturalnames=true]{hyperref}

\title{On {R}eedy model categories}
\author{Clark Barwick}
\address{Matematisk Institutt\\
Universitetet i Oslo\\
Boks 1053 Blindern\\
0316 Oslo\\
Norge}
\email{clarkbar@gmail.com}
\date{17 August 2007}

\newcommand{\fromto}[2]{\xymatrix@1@C=18pt{{#1}\ar[r]&{#2}}}
\newcommand{\goesto}[2]{\xymatrix@1@C=12pt{{#1}\,\ar@{|->}[r]&{#2}}}
\newcommand{\adjunct}[4]{{#1}:\xymatrix@1@C=18pt{{#2}\ar@<0.5ex>[r]&{#3}\ar@<0.5ex>[l]}:{#4}}

\begin{document}

\maketitle
\thispagestyle{empty}

The sole purpose of this note is to introduce some elementary results on the structure and functoriality of Reedy model categories. Presumably experts will have known most of the results produced here for some time, but it may be the case that there are one or two results that have not become part of the conventional wisdom.

This very brief note was culled --- with only mild changes --- from my forthcoming books \cite{barwick_book1} and \cite{barwick_book2} on higher categories and weak enrichments. I have decided to make some results available separately, in deference those who apparently wish to use the some of the techniques before the long process of editing the books is complete.

After a very brief reprise of well-known facts about the Reedy model structure, I give a very useful little criterion to determine whether composition with a morphism of Reedy categories determines a left or right Quillen functor. I then give three easy inheritance results, and the paper concludes with a somewhat more difficult inheritance result, providing conditions under which the Reedy model structure on diagrams valued in a symmetric monoidal model category is itself symmetric monoidal.

Thanks to J. Bergner, P. A. {\O}stv{\ae}r, and B. Toën for persistent encouragement and hours of interesting discussion. Thanks especially to M. Spitzweck for a profound and lasting impact on my work; were it not for his insights and questions, there would be nothing for me to report here or anywhere else.

\setcounter{tocdepth}{2}
\tableofcontents

\section{Inverse, direct, and Reedy categories} Suppose $\XX$ a universe, $\MM$ a model $\XX$-category. The Reedy model structure on the category of $s\MM$ of simplicial objects of $\MM$ is well-known in the context of resolutions, but in fact the Reedy model structure for categories of diagrams indexed by \emph{any} Reedy category has significant applications in homotopy coherent algebra as well.

I begin by reviewing some definitions and results of \cite[\S 5.1]{MR99h:55031}. Let $\XX$ be a universe.

\begin{dfn}\label{dfn:directinverse} Suppose $A$ an $\XX$-small category, $\lambda$ an $\XX$-small ordinal.
\begin{enumerate}[(\ref{dfn:directinverse}.1)]
\item For any $\XX$-small category $A$, a functor $d:\fromto{A}{\lambda}$ is called a \emph{linear extension} of $A$ if it refects identities, that is, if a morphism $f$ of $A$ is an identity if and only if $d(f)$ is.
\item An $\XX$-small category $A$ is said to be a \emph{direct category} if there exists a linear extension $d:\fromto{A}{\lambda}$.
\item An $\XX$-small category $A$ is said to be an \emph{inverse category} if $A^{\op}$ is a direct category.
\end{enumerate}
\end{dfn}

\begin{thm}\label{thm:projdirectinjinverse} Suppose $\MM$ a model $\XX$-category.
\begin{enumerate}[(\ref{thm:projdirectinjinverse}.1)]
\item For any direct category $A$, the functor category $\MM^A$ has its projective model structure, in which the weak equivalences and fibrations are defined objectwise.
\item For any inverse category $A$, the functor category $\MM^A$ has its injective model structure, in which the weak equivalences and cofibrations are defined objectwise.
\end{enumerate}
\begin{proof} This is \cite[Theorem 5.1.3]{MR99h:55031}.
\end{proof}
\end{thm}

\begin{prp}\label{prp:leftrightKan} Suppose $f:\fromto{A}{B}$ a functor of $\XX$-small categories, $\MM$ a model category.
\begin{enumerate}[(\ref{prp:leftrightKan}.1)]
\item If $A$ and $B$ are direct categories, then the adjunction
\begin{equation*}
\adjunct{f_!}{\MM^A}{\MM^B}{f^{\star}}
\end{equation*}
is a Quillen adjunction between the projective model categories.
\item If $A$ and $B$ are inverse categories, then the adjunction
\begin{equation*}
\adjunct{f^{\star}}{\MM^B}{\MM^A}{f_{\star}}
\end{equation*}
is a Quillen adjunction between the injective model categories.
\end{enumerate}
\begin{proof} It is obvious that $f^{\star}$ preserves any types of morphisms that are defined objectwise.
\end{proof}
\end{prp}

\begin{dfn}\label{dfn:latchmatch} Suppose $C$ any $\XX$-complete and $\XX$-cocomplete $\XX$-category.
\begin{enumerate}[(\ref{dfn:latchmatch}.1)]
\item\label{item:latch} Suppose $A$ a direct category, $\alpha$ and object of $A$.
\begin{enumerate}[(\ref{dfn:latchmatch}.\ref{item:latch}.1)]
\item The \emph{latching category} at $\alpha$ is the full subcategory $\partial(A/\alpha)$ of the category $(A/\alpha)$ consisting of the nonidentity morphisms $\fromto{\beta}{\alpha}$. There are two forgetful functors:
\begin{equation*}
F_{\alpha}:\fromto{(A/\alpha)}{A}\textrm{\qquad and\qquad}\partial F_{\alpha}:\fromto{\partial(A/\alpha)}{A}.
\end{equation*}
\item The \emph{latching functor} $L_{\alpha}$ for $C$ is the composite functor
\begin{equation*}
\xymatrix@C=24pt{C^A\ar[r]^-{\partial F_{\alpha}^{\star}}&C^{\partial(A/\alpha)}\ar[r]^-{\colim}&C},
\end{equation*}
and the image of a diagram $X:\fromto{A}{C}$ is called the \emph{latching object $L_{\alpha}X$ of $X$ at $\alpha$.}
\end{enumerate}
\item\label{item:match} Suppose $A$ an inverse category, $\alpha$ and object of $A$.
\begin{enumerate}[(\ref{dfn:latchmatch}.\ref{item:match}.1)]
\item The \emph{matching category} at $\alpha$ is the opposite category $\partial(\alpha/A):=(\partial(A^{\op}/\alpha))^{\op}$ of the latching category at $\alpha$ for $A^{\op}$. There are two forgetful functors:
\begin{equation*}
F^{\alpha}:\fromto{(\alpha/A)}{A}\textrm{\qquad and\qquad}\partial F^{\alpha}:\fromto{\partial(\alpha/A)}{A}.
\end{equation*}
\item The \emph{matching functor} $M^{\alpha}$ for $C$ is the composite functor
\begin{equation*}
\xymatrix@C=24pt{C^A\ar[r]^-{\partial F^{\alpha,\star}}&C^{\partial(\alpha/A)}\ar[r]^-{\lim}&C},
\end{equation*}
and the image of a diagram $X:\fromto{A}{C}$ is called the \emph{matching object $M^{\alpha}X$ of $X$ at $\alpha$.}
\end{enumerate}
\end{enumerate}
\end{dfn}

\begin{prp}\label{prp:cofsprojfibsinj} Suppose $\MM$ a model $\XX$-category.
\begin{enumerate}[(\ref{prp:cofsprojfibsinj}.1)]
\item For any direct category $A$, a morphism $\fromto{X}{Y}$ of the functor category $\MM^A$ is a cofibration or trivial cofibration in the projective model structure if and only if for any object $\alpha$ of $A$, the induced morphism
\begin{equation*}
\fromto{X_{\alpha}\coprd^{L_{\alpha}X}L_{\alpha}Y}{Y_{\alpha}}
\end{equation*}
is so.
\item For any inverse category $A$, a morphism $\fromto{X}{Y}$ of the functor category $\MM^A$ is a fibration or trivial fibration in the projective model structure if and only if for any object $\alpha$ of $A$, the induced morphism
\begin{equation*}
\fromto{X_{\alpha}}{M^{\alpha}X\times_{M^{\alpha}Y}Y_{\alpha}}
\end{equation*}
is so.
\end{enumerate}
\begin{proof} This is \cite[Theorem 5.1.3]{MR99h:55031}.
\end{proof}
\end{prp}

\begin{dfn}\label{dfn_Reedycategory} A \emph{Reedy category} consists of the following data:
\begin{enumerate}[(\ref{dfn_Reedycategory}.A)]
\item an $\XX$-small category $A$,
\item two lluf subcategories $A^{\to}$ and $A^{\gets}$ of $A$, and
\item a functorial factorization of every morphism into a morphism of $A^{\gets}$ followed by a morphism of $A^{\to}$.
\end{enumerate}
These data are subject to the following condition: there exist an ordinal $\lambda$ and two linear extensions $\fromto{A^{\to}}{\lambda}$ and $\fromto{(A^{\gets})^{\op}}{\lambda}$ such that the diagram
\begin{equation*}
\xymatrix@C=18pt@R=12pt{
&A^{\to}\ar[dr]&\\
\Obj A\ar[ur]\ar[dr]&&\lambda\\
&(A^{\gets})^{\op}\ar[ur]&
}
\end{equation*}
commutes. Write $i^{\to}$ (respectively, $i^{\gets}$) for the inclusion $\fromto{A^{\to}}{A}$ (resp., for the inclusion $\fromto{A^{\gets}}{A}$).
\end{dfn}

\begin{nul} In other words, a Reedy category consists of a category with a degree function on its objects, so that any morphism can be factored in a functorial fashion as a morphism that decreases the degree followed by a morphism that increases the degree.
\end{nul}

\begin{lem}\label{lem:AReedyAopReedy} If $A$ is a Reedy category, then $A^{\op}$ is as well, with $(A^{\op})^{\to}:=(A^{\gets})^{\op}$ and $(A^{\op})^{\gets}:=(A^{\to})^{\op}$.
\begin{proof} The unique factorization for $A$ will work for $A^{\op}$.
\end{proof}
\end{lem}

\begin{lem}\label{lem:AmodgammaReedy} Suppose $A$ a Reedy category, $C$ an arbitrary category, and $\fromto{A}{C}$ a fully faithful functor.
\begin{enumerate}[(\ref{lem:AmodgammaReedy}.1)]
\item For any object $\gamma$ of $C$, the slice category $(A/\gamma)$ is a Reedy category, wherein $(A/\gamma)^{\to}$ (respectively, $(A/\gamma)^{\gets}$) is the lluf subcategory consisting of those morphisms mapping to $A^{\to}$ (resp., to $A^{\gets}$) under the obvious forgetful functor $\fromto{(A/\gamma)}{A}$.
\item For any object $\gamma$ of $C$, the slice category $(\gamma/A)$ is a Reedy category, wherein $(\gamma/A)^{\to}$ (respectively, $(\gamma/A)^{\gets}$) is the lluf subcategory consisting of those morphisms mapping to $A^{\to}$ (resp., to $A^{\gets}$) under the obvious forgetful functor $\fromto{(\gamma/A)}{A}$.
\end{enumerate}
\begin{proof} By the previous lemma, it suffices to show that $(A/\gamma)$ is a Reedy category. It is clear that the composites $\xymatrix@1@C=18pt{(A/\gamma)^{\to}\ar[r]&A^{\to}\ar[r]&\lambda}$ and $\xymatrix@1@C=18pt{(A/\gamma)^{\gets,\op}\ar[r]&A^{\gets,\op}\ar[r]&\lambda}$ are linear extensions. The unique functorial factorization for $A$ gives a unique functorial factorization for $(A/\gamma)$.
\end{proof}
\end{lem}

\begin{thm} Suppose $A$ an $\XX$-small Reedy category. Then for any model $\XX$-category $\MM$, the diagram category $\MM^A$ has its \emph{Reedy model structure}, in which a morphism $\phi:\fromto{X}{Y}$ is a weak equivalence, cofibration, or fibration if and only if both $i^{\to,\star}\phi$ in $M^{A^{\to}}$ and $i^{\gets,\star}\phi$ in $M^{A^{\gets}}$ are so.
\begin{proof} This is \cite[Theorem A]{reedy} and \cite[Theorems 15.3.4 and 15.3.15]{MR2003j:18018}.
\end{proof}
\end{thm}

\begin{nul} Note in particular that the weak equivalences are the objectwise weak equivalences.
\end{nul}

\begin{lem}\label{lem_QuillensM} The Reedy model structure is functorial in the model category; that is, suppose $A$ an $\XX$-small Reedy category, $\MM$ and $\NN$ model $\XX$-categories, and $F:\fromto{\MM}{\NN}$ a left Quillen functor. Then the induced functor $\fromto{\MM^A}{\NN^A}$ --- which will also be denoted $F$ --- is left Quillen as well.
\begin{proof} Since $F$ is a left adjoint, it commutes with all latching functors.
\end{proof}
\end{lem}

\section{Left and right fibrations of Reedy categories}{I now address the question of the functoriality of the Reedy model structure in the Reedy category. That is, I will describe the circumstances under which a functor $\fromto{A}{B}$ induces a Quillen adjunction between $\MM^A$ and $\MM^B$.}

\begin{dfn}\label{dfn:morphofReedycats} Suppose $A$ and $B$ $\XX$-small Reedy categories.
\begin{enumerate}[(\ref{dfn:morphofReedycats}.1)]
\item A \emph{morphism} $f:\fromto{A}{B}$ is a strictly commutative diagram of functors
\begin{equation*}
\xymatrix@C=18pt@R=24pt{
A^{\to}\ar[d]\ar[r]&B^{\to}\ar[d]\\
A\ar[r]&B\\
A^{\gets}\ar[u]\ar[r]&B^{\gets}.\ar[u]
}
\end{equation*}
\item A morphism $f:\fromto{A}{B}$ is a \emph{left fibration} if for any model $\XX$-category $\MM$, the adjunction
\begin{equation*}
\adjunct{f_!}{\MM^A}{\MM^B}{f^{\star}}
\end{equation*}
is a Quillen adjunction. If $B=\star$, then one says that $A$ is \emph{left fibrant}.
\item A morphism $f:\fromto{A}{B}$ is a \emph{right fibration} if for any model $\XX$-category $\MM$, the adjunction
\begin{equation*}
\adjunct{f^{\star}}{\MM^B}{\MM^A}{f_{\star}}
\end{equation*}
is a Quillen adjunction. If $B=\star$, then one says that $A$ is \emph{right fibrant}.
\end{enumerate}
\end{dfn}

\begin{nul}\label{nul_Reedybetweeninjproj} A Reedy model category is thus left (respectively, right) fibrant if and only if it has fibrant (resp., cofibrant) constants in the sense of Hirschhorn \cite[Definition 15.10.1]{MR2003j:18018}. The notion of a left or right fibration is merely a relative version of Hirschhorn's concepts.

The Reedy model structure lives between the injective model structure and projective model structure on $\MM^A$, if they exist. That is, the identity functor induces a right Quillen functor $\fromto{\MM^A_{\mathrm{Reedy}}}{\MM^A_{\mathrm{proj}}}$ and a left Quillen functor $\fromto{\MM^A_{\mathrm{Reedy}}}{\MM^A_{\mathrm{inj}}}$. If $A$ is direct (respectively, inverse), then the former (resp., latter) of these is an isomorphism of model categories. If $A$ is left fibrant (respectively, right fibrant), then the fact that the constant functor is right (resp., left) Quillen is an indication that the Reedy model structure is closer to the projective (resp., injective) model structure.
\end{nul}

\begin{lem}\label{lem:dirleftfibinvrightfib} If $A$ and $B$ are direct (respectively, inverse) categories, any morphism $f:\fromto{A}{B}$ is a left (resp., right) fibration.
\begin{proof} Immediate from \ref{prp:leftrightKan}.
\end{proof}
\end{lem}

\begin{lem} For any Reedy categories $A$ and $B$, a morphism $f:\fromto{A}{B}$ is a left fibration if and only if the functor $f^{\op}:\fromto{A^{\op}}{B^{\op}}$ is a right fibration.
\begin{proof} This follows from \ref{lem:AReedyAopReedy}.
\end{proof}
\end{lem}

\begin{lem} For any Reedy categories $A$ and $B$, a morphism $f:\fromto{A}{B}$ is a left (respectively, right) fibration if and only if the functor $f^{\gets}:\fromto{A^{\gets}}{B^{\gets}}$ (resp., the functor $f^{\to}:\fromto{A^{\to}}{B^{\to}}$) is so.
\begin{proof} By the previous lemma, it suffices to prove the statement for left fibrations. Since Reedy cofibrations and fibrations are defined be restriction to the direct and inverse subcategories, it follows that $f$ is a left fibration if and only if $f^{\to}$ and $f^{\gets}$ are left fibrations. But $f^{\to}$ is automatically a left fibration by \ref{lem:dirleftfibinvrightfib}.
\end{proof}
\end{lem}

\begin{lem} For any Reedy categories $A$ and $B$, a morphism $f:\fromto{A}{B}$ is a left (respectively, right) fibration if and only if for any object $\beta$ of $B$, the Reedy category $(f/\beta)$ (resp., $(\beta/f)$) is left (resp., right) fibrant.
\begin{proof} Again it suffices to prove the statement for left fibrations, and by the previous lemma, it suffices to assume that $A$ and $B$ are inverse categories. Now $f$ is a left fibration if and only if, for any model category $\MM$ and any (trivial) cofibration $\phi:\fromto{X}{Y}$ of $\MM^A$, the induced morphism $f_!\phi:\fromto{f_!X}{f_!Y}$ is a (trivial) cofibration of $\MM^B$. But (trivial) cofibrations are defined objectwise; hence this is in turn equivalent to the assertion that for any model category $\MM$, any (trivial) cofibration $\phi:\fromto{X}{Y}$ of $\MM^A$, and any object $\beta$ of $B$, the morphism
\begin{equation*}
f_!\phi_{\beta}:\fromto{(f_!X)_{\beta}=\colim_{\alpha\in(f/\beta)}X_{\alpha}}{\colim_{\alpha\in(f/\beta)}Y_{\alpha}=(f_!Y)_{\beta}}
\end{equation*}
is a (trivial) cofibration of $\MM$. This is precisely the statement that the adjunction
\begin{equation*}
\adjunct{\colim}{\MM^{(f/\beta)}}{\MM}{\const}
\end{equation*}
is a Quillen adjunction, i.e., that $(f/\beta)$ is left fibrant.
\end{proof}
\end{lem}

\begin{thm} For any Reedy categories $A$ and $B$, a morphism $f:\fromto{A}{B}$ is a left (respectively, right) fibration if and only if for any object $\alpha$ of $A$ and any morphism $\fromto{f(\alpha)}{\beta}$ (resp., $\fromto{\beta}{f(\alpha)}$) of $B$, the nerve of the category $\partial(\alpha/(f^{\gets}/\beta))$ (resp., of the category $\partial((\beta/f^{\to})/\alpha)$) is either empty or connected.
\begin{proof} This now follows from the previous lemma and Hirschhorn's necessary and sufficient condition for a Reedy category to have (co)fibrant constants \cite[Proposition 15.10.2(1) and Corollary 15.10.5]{MR2003j:18018}.
\end{proof}
\end{thm}

\begin{cor}\label{lem:forgetfulfibration} Suppose $A$ a Reedy category, $C$ an arbitrary category with all finite products (respectively, finite coproducts), $\fromto{A}{C}$ a fully faithful functor. Suppose that for any object $\gamma$ of $C$, the Reedy category $(A/\gamma)$ (resp., $(\gamma/A)$) is left (resp., right) fibrant. Then for any morphism $\fromto{\gamma}{\gamma'}$ of $C$, the forgetful functor $\fromto{(A/\gamma)}{(A/\gamma')}$ (resp., $\fromto{(\gamma'/A)}{(\gamma/A)}$) is a left (resp., right) fibration.
\begin{proof} Again it suffices to prove the assertion for left fibrations. Using the characterization of the theorem, one sees that the forgetful functor $\fromto{(A/\gamma)}{(A/\gamma')}$ is a left fibration if any only if for any object $\alpha$ of $(A/\gamma')$, the Reedy category $(A/(\alpha\times\gamma'))$ is left fibrant.
\end{proof}
\end{cor}

\section{Lemmata of inheritance}{I now reiterate some familiar but nevertheless useful facts on the subject of the Reedy model structure. In particular, it inherits many good formal properties of $\MM$ (\ref{lem_Reedyproper}, \ref{lem_Reedycombin}, \ref{lem_ReedysMMenriched}, and \ref{lem_MVenrichedsMVenriched}).

Suppose $A$ a Reedy category, $\MM$ a model $\XX$-category. For the sake of consistency with the case $A=\Delta$, I consider the Reedy model structure on the category $\MM(A)$ of functors $\fromto{A^{\op}}{\MM}$.}

\begin{lem}\label{lem_Reedyproper} If $\MM$ is left (respectively, right) proper, then so is the Reedy model category $\MM(A)$.
\begin{proof} This follows immediately from the observation that the Reedy weak equivalences, cofibrations, and fibrations are in particular objectwise weak equivalences, cofibrations, and fibrations.
\end{proof}
\end{lem}

\begin{nul} The Reedy model structure on $\MM(A)$ is frequently compatible with a natural symmetric monoidal structure, which arises from the use of objects $y_{\MM}(\alpha)$ that represent evaluation at an object $\alpha$ of $A$. The category $\MM(A)$ of $\MM$-valued presheaves $Y$ on $A$ comprise the representable $\mathrm{Set}_{\XX}(A)$-valued presheaves on $C$, whose value on an object $X$ of $\MM$ is the presheaf that assigns to any object $\alpha$ of $A$ the morphisms in $\MM(A)$ from a presheaf $y(\alpha)\boxdot X$ to $Y$. Extending this correspondence to all presheaves on $A$ in the usual fashion, one arrives at a fundamental adjunction of two variables (\ref{prp_boxdotadj}) on $\MM(A)$ with $\MM$ over $\mathrm{Set}_{\XX}(A)$.
\end{nul}

\begin{ntn}\label{ntn_sCsimplicial} Suppose $X$ an object of $\MM$, and $Y:\fromto{A^{\op}}{\MM}$ a presheaf.
\begin{enumerate}[(\ref{ntn_sCsimplicial}.1)]
\item Write $\mor_{\MM(A),\boxdot}(-,Y)$ for the right Kan extension $a_{\star}Y$ of $Y$ along the opposite $a:\fromto{A^{\op}}{\mathrm{Set}_{\XX}(A)^{\op}}$ of the Yoneda embedding.
\item The copower functor
\begin{equation*}
\xymatrix@C=10pt@R=0pt{\mathrm{Set}_{\XX}\ar[r]&\MM\\
S\ar@{|->}[r]&S\cdot X}
\end{equation*}
induces the functor $-\boxdot_{\MM(A)}X:\fromto{\mathrm{Set}_{\XX}(A)}{\MM(A)}$.
\item The object $X$ corepresents a functor $\fromto{\MM}{\mathrm{Set}_{\XX}}$ and thus induces a functor $\Mor_{\MM(A),\boxdot}^s(X,-):\fromto{\MM(A)}{\mathrm{Set}_{\XX}(A)}$.
\end{enumerate}
\end{ntn}

\begin{lem}\label{prp_boxdotend} For any presheaf $K:\fromto{A^{\op}}{\mathrm{Set}_{\XX}}$ and any presheaf $Y:\fromto{A^{\op}}{\MM}$, there is an isomorphism
\begin{equation*}
\mor_{\MM(A),\boxdot}(K,Y)\cong\int_{\alpha\in A}\mor(K_{\alpha},Y_{\alpha}).
\end{equation*}
\begin{proof} This is the usual end formula for right Kan extensions.
\end{proof}
\end{lem}

\begin{lem}\label{prp_boxdotadj} The triple $(\boxdot_{\MM(A)},\mor_{\MM(A),\boxdot},\Mor_{\MM(A),\boxdot}^s)$ is an adjunction of two variables: for any presheaf $K:\fromto{A^{\op}}{\mathrm{Set}_{\XX}}$, any object $X$ of $\MM$, and any presheaf $Y:\fromto{A^{\op}}{\MM}$, there are natural isomorphisms
\begin{eqnarray}
\Mor_{\MM}(X,\mor_{\MM(A),\boxdot}(K,Y))&\cong&\Mor_{\MM(A)}(K\boxdot_{\MM(A)}X,Y)\nonumber\\
&\cong&\Mor_{\mathrm{Set}_{\XX}(A)}(K,\Mor_{\MM(A),\boxdot}^s(X,Y)).\nonumber
\end{eqnarray}
\begin{proof} This follows from the relevant universal properties.
\end{proof}
\end{lem}

\begin{ntn} Write $y:\fromto{A}{\mathrm{Set}_{\XX}(A)}$ for the Yoneda embedding.
\end{ntn}

\begin{cor} Suppose $Y:\fromto{A^{\op}}{\MM}$ a presheaf; then for any object $\alpha$ of $A$, there is a natural isomorphism
\begin{equation*}
Y_{\alpha}\cong\mor_{\MM(A),\boxdot}(y(\alpha),Y).
\end{equation*}
\end{cor}

\begin{cor} Suppose $Y:\fromto{A^{\op}}{\MM}$ a presheaf; then for any object $\alpha$ of $A$, there is a presheaf $\partial y(\alpha):\fromto{A^{\op}}{\MM}$ and a natural isomorphism
\begin{equation*}
M^{\alpha}Y\cong\mor_{\MM(A),\boxdot}(\partial y(\alpha),Y).
\end{equation*}
\begin{proof} Set
\begin{equation*}
\partial y(\alpha):=\colim_{\alpha'\in(\alpha/A^{\gets})}y(\alpha').
\end{equation*}
Then one shows easily that $\mor_{\MM(A),\boxdot}(\partial y(\alpha),Y)$ is the desired limit.
\end{proof}
\end{cor}

\begin{ntn} For any set $K$ of morphisms of $\MM$, write
\begin{equation*}
\Lambda\Box K:=\{\fromto{(y(\alpha)\boxdot X)\sqcup^{\partial y(\alpha)\boxdot X}(\partial y(\alpha)\boxdot Y)}{y(\alpha)\boxdot Y}\ |\ \alpha\in A, [\fromto{X}{Y}]\in K\}.
\end{equation*}
\end{ntn}

\begin{lem}\label{lem_Reedycombin} For any $\XX$-small Reedy category $A$, the Reedy model category $\MM(A)$ is $\XX$-combinatorial \cite[1.3.1]{arXiv:0708.2067v1} if $\MM$ is.
\begin{proof} Since $\MM$ is $\XX$-combinatorial, it is possible to choose $\XX$-small sets of generating cofibrations and generating trivial cofibrations $I_{\MM}$ and $J_{\MM}$ such that the domains and codomains of $I_{\MM}$ (respectively, of $J_{\MM}$) are small with respect to $I_{\MM}$ (resp., to $J_{\MM}$); then $\Lambda\Box I_{\MM}$ and $\Lambda\Box J_{\MM}$ are $\XX$-small sets of generating cofibrations and generating trivial cofibrations of the Reedy model structure on $\MM(A)$ \cite[Theorem 15.6.27]{MR2003j:18018}. Local presentability is inherited by functor categories; hence $\MM(A)$ is $\XX$-combinatorial.
\end{proof}
\end{lem}

\begin{lem} For any $\XX$-small Reedy category $A$, the Reedy model category $\MM^A$ is $\XX$-tractable \cite[1.3.1]{arXiv:0708.2067v1} if $\MM$ is.
\begin{proof} I claim that if $\fromto{X}{Y}$ is a cofibration with cofibrant source in $\MM$, then $(y(\alpha)\boxdot X)\sqcup^{\partial y(\alpha)\boxdot X}(\partial y(\alpha)\boxdot Y)$ is cofibrant. Suppose that $\fromto{T}{S}$ is an objectwise trivial fibration. Then by adjunction, a morphism $\fromto{y(\alpha)\boxdot X}{S}$ has a lifting if and only if the diagram
\begin{equation*}
\xymatrix@C=18pt@R=18pt{
&T_{\alpha}\ar[d]\\
X\ar[r]&S_{\alpha}
}
\end{equation*}
has a lifiting. It follows from the cofibrancy of $X$ and \cite[Proposition 15.3.11]{MR2003j:18018} that $y(\alpha)\boxdot X$ is cofibrant in the Reedy model structure on $\MM(A)$. It is easy to see by a similar argument that $\fromto{\partial y(\alpha)\boxdot X}{\partial y(\alpha)\boxdot Y}$ is a cofibration, so $\fromto{y(\alpha)\boxdot X}{(y(\alpha)\boxdot X)\sqcup^{\partial y(\alpha)\boxdot X}(\partial y(\alpha)\boxdot Y)}$ is a cofibration, whence follows the claim, and thus the lemma.
\end{proof}
\end{lem}

\section{Reedy diagrams in a symmetric monoidal model category}{Suppose now $A$ an $\XX$-small Reedy category and $(\MM,\otimes_{\MM},\MOR_{\MM}^{\MM})$ a symmetric monoidal model $\XX$-category.}

\begin{ntn}\label{ntn_simplsM} Suppose $X$ and $Y$ objects of $\MM(A)$ and $Z$ an object of $\MM$. Set
\begin{eqnarray}
\MOR_{\MM(A)}^{\MM}(X,Y)&:=&\int_{\alpha\in A^{\op}}\MOR_{\MM}(X_{\alpha},Y_{\alpha}),\nonumber\\
(Z\otimes_{\MM(A)}^{\MM}X)_{\alpha}&:=&Z\otimes_{\MM}X_{\alpha},\nonumber\\
\mor_{\MM(A)}^{\MM}(Z,Y)_{\alpha}&:=&\MOR_{\MM}(Z,Y_{\alpha}),\nonumber
\end{eqnarray}
for any object $\alpha$ of $A$. This gives $\MM(A)$ the structure of an $\MM$-category.
\end{ntn}

\begin{lem}\label{lem_ReedysMMenriched} With the $\MM$-structure of \ref{ntn_simplsM}, the Reedy model category $\MM(A)$ is an $\MM$-model category.
\begin{proof} To verify the pushout-product axiom, suppose $f:\fromto{Z}{Z'}$ a cofibration of $\MM$, and $i:\fromto{X}{Y}$ a cofibration of $\MM(A)$; then for any object $\alpha$ of $A$, the morphism
\begin{equation*}
\xymatrix@R=18pt{((Z\otimes Y)\sqcup^{Z\otimes X}(Z'\otimes X))_{\alpha}\sqcup^{M^{\alpha}((Z\otimes Y)\sqcup^{Z\otimes X}(Z'\otimes X))}M^{\alpha}(Z'\otimes Y)\ar[d]\\
(Z'\otimes Y)_{\alpha}}
\end{equation*}
is isomorphic to the morphism
\begin{equation*}
\fromto{(Z\otimes Y_{\alpha})\sqcup^{Z\otimes(X_{\alpha}\sqcup^{M^{\alpha}X}M^{\alpha}Y)}(Z'\otimes(X_{\alpha}\sqcup^{M^{\alpha}X}M^{\alpha}Y)}{Z'\otimes Y_{\alpha}},
\end{equation*}
which, by the pushout-product axiom for $\MM$, is a cofibration that is trivial if either $f$ or $i$ is.
\end{proof}
\end{lem}

\begin{cor}\label{lem_MVenrichedsMVenriched} If, in addition, $\MM$ is a model $\VV$-category \cite[3.2.4]{arXiv:0708.2067v1} for some symmetric monoidal model $\XX$-category $\VV$, then $\MM(A)$ is also.
\end{cor}

\begin{lem}\label{lem_Deltacof} There is a functor $y_{\MM}:\fromto{A}{\MM(A)}$ such that for any object $\alpha$ of $A$ and any $Y:\fromto{A^{\op}}{\MM}$, there is a canonical isomorphism
\begin{equation*}
\MOR_{\MM(A)}(y_{\MM}(\alpha),Y)\cong Y_{\alpha}.
\end{equation*}
Moreover, if the unit $1_{\MM}$ for the symmetric monoidal structure on $\MM$ is cofibrant, then for every such object $\alpha$, $y_{\MM}(\alpha)$ is cofibrant.
\begin{proof} Set $y_{\MM}(\alpha):=y(\alpha)\boxdot 1_{\MM}$. The first part of the result now follows from the enriched Yoneda lemma.
\end{proof}
\end{lem}

\begin{cor}\label{cor_QuillenDeltacof} If $F:\fromto{\MM}{\NN}$ is a left Quillen functor of symmetric monoidal model $\XX$-categories such that $F(1_{\MM})\cong 1_{\NN}$,\footnote{Note that one need \emph{not} assume that $F$ itself is symmetric monoidal.} then $F(y_{\MM}(\alpha))\cong y_{\MM}(\alpha)$ for every object $\alpha$ of $A$.
\end{cor}

\begin{cor} For any object $\alpha$ of $A$, there is a simplicial object $\partial y_{\MM}(\alpha)$ of $\MM$ --- which is cofibrant if $1_{\MM}$ is --- such that for any $Y:\fromto{A^{\op}}{\MM}$, there is a canonical isomorphism
\begin{equation*}
\MOR_{\MM(A)}(\partial y_{\MM}(\alpha),Y)\cong M^{\alpha}Y.
\end{equation*}
\end{cor}

\begin{nul} The exterior tensor product $\fromto{\MM(A)\times\MM(A)}{\MM(A\times A)}$ is part of a Quillen adjunction of two variables. In order to see this, I quote the following result of Hirschhorn.
\end{nul}

\begin{thm}[Hirschhorn, \protect{\cite[Theorem 15.5.2]{MR2003j:18018}}]\label{thm:ReedyAxAReedyReedy} The category $A\times A$ has a natural Reedy category structure, for which the Reedy model structure on $\MM(A\times A)$ coincides with the ``Reedy-Reedy'' model structure on $\MM(A)(A)$.
\end{thm}

\begin{ntn}\label{ntn:MAMAMAA} Denote by
\begin{equation*}
\xymatrix@C=10pt@R=0pt{
\boxast_{\MM(A)}:\MM(A)\times\MM(A)\ar[r]&\MM(A\times A)\\
\MOR_{\boxast,\MM(A)}:\MM(A)^{\op}\times\MM(A\times A)\ar[r]&\MM(A)\\
\mor_{\boxast,\MM(A)}:\MM(A)^{\op}\times\MM(A\times A)\ar[r]&\MM(A)
}
\end{equation*}
the functors defined by the formul\ae
\begin{eqnarray}
(X\boxast_{\MM(A)}Y)_{(\alpha,\alpha')}&:=&X_{\alpha}\otimes_{\MM}Y_{\alpha'},\\
\MOR_{\boxast,\MM(A)}(Y,F)_{\alpha}&:=&\MOR_{\MM(A\times A)}^{\MM}((y_{\MM}(\alpha)\boxast_{\MM(A)}Y),F),\\
\mor_{\boxast,\MM(A)}(X,F)_{\alpha}&:=&\MOR_{\MM(A\times A)}^{\MM}((X\boxast_{\MM(A)}y_{\MM}(\alpha)),F),
\end{eqnarray}
for any objects $X$ and $Y$ of $\MM(A)$, any $F:\fromto{A^{\op}\times A^{\op}}{\MM}$, and any objects $\alpha,\alpha'$ of $A$.
\end{ntn}

\begin{prp} The triple $(\boxast,\MOR_{\boxast,\MM(A)},\mor_{\boxast,\MM(A)})$ is an adjunction of two variables from $\MM(A)\times\MM(A)$ to $\MM(A\times A)$.
\begin{proof} This is an easy consequence of the Fubini theorem for ends and the representability properties of $y_{\MM}(\alpha)$.
\end{proof}
\end{prp}

\begin{lem} For any pair of objects $\alpha$ and $\beta$ of $A$, there is a canonical isomorphism
\begin{equation*}
y_{\MM}(\alpha)\boxast y_{\MM}(\beta)\cong y_{\MM}(\alpha,\beta)
\end{equation*}
in the category $\MM(A\times A)$.
\begin{proof} This follows from the observation that for any $\XX$-small sets $S$ and $T$, there is a canonical isomorphism
\begin{equation*}
(S\cdot\mathbf{1})\otimes (T\cdot\mathbf{1})\cong(S\times T)\cdot\mathbf{1}
\end{equation*}
in $\MM$.
\end{proof}
\end{lem}

\begin{cor} For any object $\alpha$ of $A$ and any $F:\fromto{A^{\op}\times A^{\op}}{\MM}$, there is a canonical isomorphism
\begin{equation*}
\MOR_{\boxast,\MM(A)}(y(\alpha),F)\cong F(\alpha,-)
\end{equation*}
in $\MM(A)$.
\end{cor}

\begin{cor} For any object $\beta$ of $A$, any object $X$ of $\MM(A)$, and any $F:\fromto{A^{\op}\times A^{\op}}{\MM}$, there is a canonical isomorphism
\begin{equation*}
\MOR_{\boxast,\MM(A)}(X,F)_{\beta}\cong\MOR_{\MM(A)}^{\MM}(X,F(-,\beta))
\end{equation*}
in $\MM$.
\end{cor}

\begin{cor} For any object $\beta$ of $A$, any object $X$ of $\MM(A)$, and any $F:\fromto{A^{\op}\times A^{\op}}{\MM}$, there is a canonical isomorphism
\begin{equation*}
M_{\beta}\MOR_{\boxast,\MM(A)}(X,F)\cong\MOR_{\MM(A)}^{\MM}(X,M_{(-,\beta)}F)
\end{equation*}
in $\MM$.
\end{cor}

\begin{prp} The adjunction of two variables $(\boxast,\MOR_{\boxast,\MM(A)},\mor_{\boxast,\MM(A)})$ is a Quillen adjunction of two variables.
\begin{proof} Now suppose $i:\fromto{X}{Y}$ a cofibration of $\MM(A)$, and suppose $p:\fromto{F}{G}$ a (trivial) fibration of $\MM(A)$. Now by \ref{thm:ReedyAxAReedyReedy}, for any object $\beta$ of $A$, $p$ induces a (trivial) fibration
\begin{equation*}
\fromto{F(-,\beta)}{M_{(-,\beta)}(p)}
\end{equation*}
of $\MM(A)$, where
\begin{equation*}
M_{(-,\beta)}(p):=M_{(-,\beta)}F\times_{M_{(-,\beta)}G}G(-,\beta).
\end{equation*}
Since $\MM(A)$ is a model $\MM$-category, it follows that the induced morphism
\begin{equation*}
\xymatrix@R=18pt{
\MOR_{\MM(A)}(Y,F(-,\beta)\ar[d]\\
(\MOR_{\MM(A)}(Y,M_{(-,\beta)}(p))\times_{\MOR_{\MM(A)}(X,M_{(-,\beta)}(p))}\MOR_{\MM(A)}(X,F(-,\beta))
}
\end{equation*}
is a fibration of $\MM$, which is trivial if either $i$ or $p$ is. But this morphism is isomorphic to the morphism $\fromto{U_{\beta}}{M_{\beta}U\times_{M_{\beta}V}V_{\beta}}$, wherein
\begin{eqnarray}
U&:=&\MOR_{\boxast,\MM(A)}(Y,F);\nonumber\\
V&:=&\MOR_{\boxast,\MM(A)}(X,F)\times_{\MOR_{\boxast,\MM(A)}(X,G)}\MOR_{\boxast,\MM(A)}(Y,G).\nonumber
\end{eqnarray}
Hence $\fromto{U}{V}$ is a fibration of $\MM(A)$, which is trivial if $i$ or $p$ is.
\end{proof}
\end{prp}

\begin{prp} Suppose
\begin{equation*}
\adjunct{F}{\MM(A\times A)}{\MM(A)}{U}
\end{equation*}
a Quillen adjunction. Then the triple $(\otimes_{\MM(A),F},\MOR_{\MM(A),F},\mor_{\MM(A),F})$, defined by the formul\ae
\begin{eqnarray}
X\otimes_{\MM(A),F}Y&:=&F(X\boxast_{\MM(A)}Y),\\
\MOR_{\MM(A),F}(Y,Z)&:=&\MOR_{\boxast,\MM(A)}(Y,UZ),\\
\mor_{\MM(A),F}(X,Z)_{\alpha}&:=&\mor_{\boxast,\MM(A)}(X,UZ),
\end{eqnarray}
for any objects $X$, $Y$, and $Y$ of $\MM(A)$, is a Quillen adjunction of two variables from $\MM(A)\times\MM(A)$ to $\MM(A)$.
\begin{proof} Suppose $i:\fromto{U}{V}$ and $j:\fromto{X}{Y}$ cofibrations of $\MM(A)$. Then $i\Box_{\MM(A),F}j=F(i\Box_{\boxast,\MM(A)}j)$, which is a cofibration that is trivial if either $i$ or $j$ is.
\end{proof}
\end{prp}

\begin{cor} Suppose $A$ monoidal, with a structure
\begin{equation*}
\circ:\fromto{A\times A}{A}
\end{equation*}
that defines a right fibration of Reedy categories. Then the Day convolution product
\begin{equation*}
\otimes_{\MM(A),\circ}:=\circ_!(-\boxast_{\MM(A)}-)
\end{equation*}
is part of a Quillen adjunction of two variables on $\MM(A)$.
\end{cor}

\begin{thm}\label{thm:symmonReedysymmon} If $A$ is left fibrant and the morphisms of $A^{\gets}$ are epimorphisms, then the \emph{diagonal} symmetric monoidal structure given by
\begin{eqnarray}
(X\otimes_{\MM(A)}^{\MM(A)}Y)_{\alpha}&:=&X_{\alpha}\otimes Y_{\alpha},\\
\MOR_{\MM(A)}^{\MM(A)}(X,Y)_{\alpha}&:=&\MOR_{\MM(A)}^{\MM}(y_{\MM}(\alpha)\otimes_{\MM(A)}^{\MM(A)}X,Y),
\end{eqnarray}
for any objects $X$ and $Y$ of $\MM(A)$ and an object $\alpha$ of $A$, gives $\MM(A)$ the structure of a symmetric monoidal model category.
\begin{proof} The unit axiom follows from the fact that the constant functor is symmetric monoidal and preserves cofibrant objects and equivalences.

It now suffices to show that the diagonal functor $\Delta:\fromto{A}{A\times A}$ is a left fibration. This is equivalent to showing, for any objects $\alpha,\beta,\gamma$ of $A$, and any pair of morphisms $\fromto{\alpha}{\beta}$ and $\fromto{\alpha}{\gamma}$ of $A^{\gets}$, that the nerve of the category $\partial(\alpha/(\Delta^{\gets}/(\beta,\gamma))/)$ is either empty or connected. Since $A$ is left fibrant, the nerve of the category $\delta(\alpha/A^{\gets})$ is either empty or connected. Hence if
\begin{equation*}
\xymatrix@C=10pt@R=10pt{
&(\delta,\delta)\ar[dr]&\\
(\alpha,\alpha)\ar[ur]\ar[rr]\ar[dr]&&(\beta,\gamma)\\
&(\epsilon,\epsilon)\ar[ur]&
}
\end{equation*}
is a commutative diagram of $A^{\gets}\times A^{\gets}$, then there exists a zig-zag of morphisms of $A^{\gets}\times A^{\gets}$ connecting $(\delta,\delta)$ to $(\epsilon,\epsilon)$ \emph{under} $(\alpha,\alpha)$. To see that these morphisms are morphisms over $(\beta,\gamma)$ as well, we can, without loss of generality, suppose that there is a morphism $\fromto{(\delta,\delta)}{(\epsilon,\epsilon)}$ of $((\alpha,\alpha)/(A^{\gets}\times A^{\gets}))$. Hence the left half and the exterior square of the diagram
\begin{equation*}
\xymatrix@C=10pt@R=10pt{
&(\delta,\delta)\ar[dd]\ar[dr]&\\
(\alpha,\alpha)\ar[ur]\ar[dr]&&(\beta,\gamma)\\
&(\epsilon,\epsilon)\ar[ur]&
}
\end{equation*}
commute. But since $\fromto{(\alpha,\alpha)}{(\epsilon,\epsilon)}$ is an epimorphism, it follows that the left half of this diagram commutes as well.
\end{proof}
\end{thm}

\begin{exm} For any $\XX$-small simplicial set or category $K$, the category $\MM(\Delta/K)$ is symmetric monoidal with the diagonal symmetric monoidal structure.
\end{exm}

\bibliographystyle{../../smfplain}
\bibliography{../../math,../../egasga}

\end{document}